# Multiple-Layer Parking with Screening

Sjoert Fleurke and Aernout C. D. van Enter

*Abstract*—In this article a multilayer parking system with screening of size $n = 3$ is studied with a focus on the time-dependent particle density. We prove that the asymptotic limit of the particle density increases from an average density of 1/3 on the first layer to the value of $(10 - \sqrt{5})/19 \approx 0.4086$ in higher layers.

## I. Introduction

SUPPOSE we have a lattice $\mathcal{L}(x,r)$ consisting of sites $(x,r)$ with positions $x \in \{-2,-1,0,1,2\}$ and heights $r \in \mathbb{N}$. At each position particles arrive according to independent Poisson processes $N_t(x)$. We impose boundary conditions $N_t(-2) = N_t(2) = 0$. The particles pile up across the layers but they are not allowed to "interfere" with particles earlier deposited in neighboring sites at the same layer. In other words, the horizontal distance between two particles has to be at least 2. Furthermore, in this model the particles are not allowed to pass earlier deposited particles. As a consequence a new particle is always deposited in the layer above the highest layer that rejected it. This model property is sometimes called "screening" (see Fig. 1).

Our model can be formulated more precisely in the following way.

1) The state-space is $\mathcal{F} \coloneqq (\mathcal{L}, \mathbb{N}^+)^{\{0,1\}}$.

2) The process $\kappa_t(x,r) = 1$ if there is a particle at $(x,r)$ at time t and 0 otherwise.

3) When a particle arrives at site x at time t, it will be deposited at $h_t(x) \coloneqq 1 + \max\{r : \exists_{y \in N_x}, \kappa_t(y,r) = 1\}$, where neighborhood set $N_x$ consists of site x and the sites with distance 1 from it.

The density $\rho_t(x,r)$ of a site at $(x,r) \in \mathcal{L}$ is defined as the expectation of the occupancy of that site at time $t$, i.e. $\rho_t(x,r) = E\kappa_t(x,r)$. The end-density of a site is $\rho_\infty(x,r)$.

The majority of the existing literature in which discrete

This work was supported by the Radio Communications Agency, Netherlands.
S. R. Fleurke is with the Radiocommunications Agency Netherlands, Postbus 450, 9700 AL, Groningen, The Netherlands. (corresponding author: sjoert.fleurke@agentschaptelecom.nl).
A. C. D. van Enter is with the Johann Bernoulli Institute for Mathematics and Computer Science, University of Groningen, Nijenborgh 9, 9747 AG, Groningen, The Netherlands. (e-mail: a.c.d.van.enter@rug.nl).

parking is analytically treated is about monolayer models [1, 2, 3], while most literature about multi-layer models is based on simulations [4, 5]. However, recently there is some interest in analytical results on multilayer parking models. In [6], for example, it was shown that in an infinite parking system the second layer has a higher end-density than the first layer for models both with and without screening. Analytical formulas for the time-dependent densities were derived for small parking systems without screening in [7] and [8]. It was found that the end-density in the case of a system of size three tends to exactly $\frac{1}{2}$ for high layers. It is conjectured that the same counts for bigger parking systems. In [9] density formulas are calculated for the model with screening in the case of infinite-sized regular and random trees. Contrary to the model without screening the layer densities turn out to decrease with the layer number. In [10] it was proven that the end-density of an infinite parking system with screening tends to a value that lies between $\frac{1}{k^*}$ and $\frac{1}{2}$, where $k^*$ is such that $\left(2\frac{e}{k^*}\right)^{k^*} - e = 0$, which means that $0.232 < \rho < 0.500$. To our knowledge a precise value for the end-densities of a system with screening of any size is yet to be found.

In this paper we continue the work on calculating the particle densities in a small multi-layer parking model. We hope our result will lead to further insights also in systems with bigger sizes and systems with neighborhoods of cardinality greater than 2.

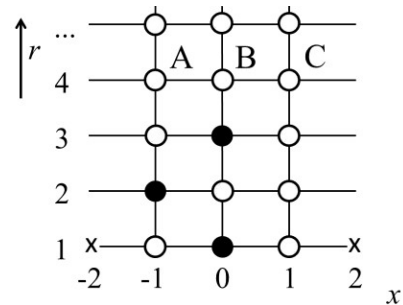

Fig. 1: Parking on a lattice with 3 positions where parking is allowed. In this picture 3 particles have arrived at positions $x_1 = 0, x_2 = -1$, and $x_3 = 0$. The next particle will be deposited at position A, B, or C depending on its x-position. Note that in the case of the parking model without screening a particle at position 1 would not be deposited at C(1,4) but at position (x, r) = (1, 2). The '×' symbols at positions -2 and 2 are indicating that at those positions no particles are dropped during the process.

## II. TIME-DEPENDENT PARTICLE DENSITIES

This section is dedicated to the calculation of the time-dependent density of vertex $(0, r)$ at time $t$. The result is displayed in the following theorem.

*Theorem 1:* Consider parking with screening on a lattice $\mathcal{L}(x, r)$ consisting of sites $(x, r)$ with positions $x \in \{-2, -1, 0, 1, 2\}$ and heights $r \in \mathbb{N}$ and boundary conditions $N_t(-2) = N_t(2) = 0$. The time-dependent particle density at layer $r$ is given by:

$$\rho_t(0, r) = \sum_{i=1}^{r} \sum_{d_1+d_2+\cdots+d_i=r} \prod_{j=1}^{i} \left[ \frac{2}{3} \sum_{k=0}^{d_j-2} \binom{d_j+k-1}{k} \left(\frac{1}{3}\right)^{d_j+k-1} + \binom{2d_j-2}{d_j-1} \left(\frac{1}{3}\right)^{2d_j-1} \right] \left( 1 - e^{-t} \sum_{l=0}^{i} \frac{t^l}{l!} \right) \quad (1)$$

Using Theorem 1 we find the following densities for the first 4 layers:

$$\rho_t(0, 1) = \frac{1}{3} - \left(\frac{1}{3} + \frac{1}{3}t\right) e^{-t} \quad (2)$$

$$\rho_t(0, 2) = \frac{11}{27} - \left(\frac{11}{27} + \frac{11}{27}t + \frac{1}{9}\frac{t^2}{2}\right) e^{-t} \quad (3)$$

$$\rho_t(0, 3) = \frac{11}{27} - \left(\frac{11}{27} + \frac{11}{27}t + \frac{19}{81}\frac{t^2}{2!} + \frac{1}{27}\frac{t^3}{3!}\right) e^{-t} \quad (4)$$

$$\rho_t(0, 4) = \frac{893}{2187} - \left(\frac{893}{2187} + \frac{893}{2187}t + \frac{229}{729}\frac{t^2}{2!} + \frac{1}{27}\frac{t^3}{3!} + \frac{1}{81}\frac{t^4}{4!}\right) e^{-t} \quad (5)$$

Theorem 1 can be proven using the fact that the process has a renewal structure. Every particle arrival at the center vertex $x = 0$ counts as a renewal. Between every arrival at the center there have been zero or more arrivals at the neighboring sites. The vertical distance between two consecutive particles arriving in the center is thus determined by the maximum of the number of arriving particles at the neighboring sites. More precisely, if $t_n$ denotes the arrival time of the $n^{th}$ particle at the center and $N_{t_1, t_2}(x)$ denotes the number of arrivals at $x$ between time $t_1$ and $t_2$, then the vertical distance between the $n^{th}$ and the $(n+1)^{st}$ consecutively arrived particles at $x = 0$ is equal to

$$\psi_n := 1 + \max\{N_{t_{n-1}, t_n}(-1), N_{t_{n-1}, t_n}(1)\}. \quad (6)$$

This means that, for example, the probability that the first center particle is deposited at height $r = 5$, is equal to $P(\psi_1 = 5) = P(\max\{N_{0, t_1}(-1), N_{0, t_1}(1)\} = 4)$.

More generally we can write:

$$\rho_t(0, r) = P\left( \sum_{j=1}^{N_t(0)} \psi_j = r \right) \quad (7)$$

This leads to:

$$\rho_t(0, r) = \int_0^t \sum_{i=1}^{r} P\left( \sum_{j=1}^{N_t(0)} \psi_j = r \middle| N_t(0) = i \right) \cdot P(N_t(0) = i) du \quad (8)$$

Note that the $(\psi_n)_{n \in \mathbb{N}^+}$ are independently and identically distributed. Also remember that $N_t(0)$ is Poisson distributed, so that (4) can be rewritten as

$$\rho_t(0, r) = \int_0^t \sum_{i=1}^{r} \sum_{d_1+d_2+\cdots+d_i=r} P\left( \bigcap_{j=1}^{i} \psi_j = d_j \right) \cdot e^{-u} \frac{u^i}{i!} du \quad (9)$$

$$= \int_0^t \sum_{i=1}^{r} \sum_{d_1+d_2+\cdots+d_i=r} \prod_{j=1}^{i} P(\psi_j = d_j) e^{-u} \frac{u^i}{i!} du \quad (10)$$

$$= \sum_{i=1}^{r} \sum_{d_1+d_2+\cdots+d_i=r} \prod_{j=1}^{i} P(\psi_j = d_j) \int_0^t e^{-u} \frac{u^i}{i!} du \quad (11)$$

$$= \sum_{i=1}^{r} \sum_{d_1+d_2+\cdots+d_i=r} \prod_{j=1}^{i} P(\psi_j = d_j) \cdot \left( 1 - e^{-t} \sum_{k=0}^{i} \frac{t^k}{k!} \right) \quad (12)$$

To complete the result we need to calculate the distribution of the stochastic variable $\psi_n$.

*Lemma 1:* The distribution of $\psi_n$ is given by:

$$P(\psi_n = d) = \frac{2}{3} \sum_{k=0}^{d-2} \binom{d+k-1}{k} \left(\frac{1}{3}\right)^{d+k-1} \quad (13)$$

$$+\binom{2d-2}{d-1}\left(\frac{1}{3}\right)^{2d-1}$$

*Proof:* For this proof we use expression (6). We can calculate the distribution of $S_T = \max\{N_T(-1), N_T(1)\}$, i.e. the maximum of the number of particles that arrive at position $x = -1$ and $x = 1$ in a period of time of length $T$.

$$P(S_T = d) =$$
$$= \int_0^\infty P(\max\{N_T(-1), N_T(1)\} = d | T = u) \cdot P(T = u) du$$
$$= \int_0^\infty P(N_T(-1) = d \cap N_T(1) < d | T = u) e^{-u} du$$
$$+ \int_0^\infty P(N_T(-1) < d \cap N_T(1) = d | T = u) e^{-u} du$$
$$+ \int_0^\infty P(N_T(-1) = N_T(1) = d | T = u) e^{-u} du$$
$$= \int_0^\infty 2P(N_T(-1) = d \cap N_T(1) < d | T = u) e^{-u} du$$
$$+ \int_0^\infty [P(N_T(1) = d | T = u)]^2 e^{-u} du \quad (14)$$
$$= 2\int_0^\infty e^{-u} \frac{u^d}{d!} \sum_{k=0}^{d-1} e^{-u} \frac{u^d}{k!} e^{-u} du$$
$$+ \int_0^\infty \left[e^{-u} \frac{u^d}{d!}\right]^2 e^{-u} du$$
$$= 2\sum_{k=0}^{x-1} \frac{(d+k)!}{d!\,k!} \left(\frac{1}{3}\right)^{d+k+1}$$
$$\cdot \int_0^\infty \frac{3^{2d+1}}{(d+k)!} u^{(d+k+1)-1} e^{-3u} du$$
$$+ \frac{(2d)!}{d!\,d!}\left(\frac{1}{3}\right)^{2d+1} \int_0^\infty \frac{3^{2d+1}}{(2d)!} u^{(2d+1)-1} e^{-3u} du$$
$$= \frac{2}{3}\sum_{k=0}^{d-1} \binom{d+k}{k}\left(\frac{1}{3}\right)^{d+k} + \binom{2d}{d}\left(\frac{1}{3}\right)^{2d+1}$$

Now we take $\psi_n = 1 + S_T$ to establish our result of Lemma 1. Finally, combining Lemma 1 with equation (12) yields the formula of Theorem 1.

In Figure 2 developments in time of particle densities are displayed for several layers. It is interesting to see that in higher layers the probability of a particle hit at the center site tends to be (slightly) bigger than in lower layers.
The opposite phenomenon was observed in larger systems in [9]. In an infinite-sized system the density decreases when parking with screening is conducted.

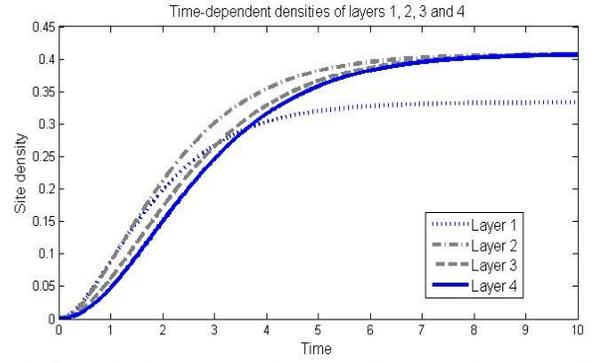

Fig. 2: Particle densities at the sites (0,r) as a function of time for r is 1, 2, 3 and 4 according to Theorem 1 and formulas (2), (3), (4) and (5) in particular.

### III. LAYER-DEPENDENT END-DENSITIES

With the above result it is straightforward to find the end-density for any layer. Take formula (1) and then let $t \to \infty$. This immediately yields:

$$\rho_\infty(0,r) = \sum_{i=1}^r \sum_{d_1+d_2+\cdots+d_i=r}$$
$$\prod_{j=1}^i \left[\frac{2}{3}\sum_{k=0}^{d_j-2}\binom{d_j+k-1}{k}\left(\frac{1}{3}\right)^{d_j+k-1} + \binom{2d_j-2}{d_j-1}\left(\frac{1}{3}\right)^{2d_j-1}\right] \quad (15)$$

In Figure 3 the end-densities are displayed for the first 4 layers based on formula (15) which yields $\frac{1}{3}, \frac{11}{27}, \frac{11}{27}$, and $\frac{893}{2187}$ respectively. These numbers have been confirmed by simulations.

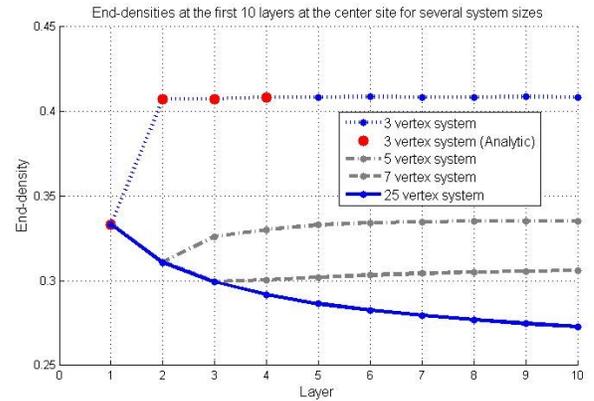

Fig. 3: End-densities as a function of the layer number based on formula (15) and simulation results. For some other bigger systems the end-densities are plotted as well. It can be seen that whereas in the 3-vertex the end-densities grow with the layer number this effect disappears in bigger systems. It is conjectured that for an infinite-sized system the end-densities of the center site decrease monotonically with the layer number.

## IV. END-DENSITY FOR HIGH LAYERS

As can be seen in Figure 3 the end-density grows with the layer number. A natural and interesting question is what the exact limiting value for this end-density could be. The purpose of this section is to prove the following theorem.

*Theorem 2.* The average end-density of the center vertex ultimately tends to the value

$$p_0 := \lim_{r \to \infty} \rho_\infty(0, r) = \frac{10 - \sqrt{5}}{19} \approx 0.408628 \quad (16)$$

The same holds for the neighboring sites: $p_{-1} = p_1 = p_0 = (10 - \sqrt{5})/19$.

One may attempt to obtain this limit from formula (15) by taking $r \to \infty$ but this seems rather difficult. Fortunately, there is an alternative, and perhaps more interesting, approach which starts with the observation that the average end-density is the inverse of 1 plus the average number of vertical consecutive empty positions. If, for example, the average run of empty sites between two occupied sites were 4, it would follow that the average occupancy is $1/(4+1) = 0.20$. In other words, the following formula holds:

$$\lim_{r \to \infty} \rho_\infty(0, r) = \frac{1}{1 + EX} \quad (17)$$

Our efforts are now concentrated on calculating the value of $EX$. We will prove the following lemma.

*Lemma 2.* The expected run of empty center sites for high layers is given by:

$$EX = 1 + \frac{1}{5}\sqrt{5} \quad (18)$$

*Proof:* If the number of particles arrived at the border sites is known, we can calculate the expected run length $EX$ as follows:

$$E(X|N = n) = \sum_{i=0}^{n} \max(i, n-i) \binom{n}{i} \left(\frac{1}{2}\right)^n \quad (19)$$

If n is an even number this can be written as:

$$\sum_{i=0}^{n/2} (n-i) \binom{n}{i} \left(\frac{1}{2}\right)^n + \sum_{i=\frac{n}{2}+1}^{n} i \binom{n}{i} \left(\frac{1}{2}\right)^n \quad (20)$$

while if n is an odd number it can be written as:

$$\sum_{i=0}^{(n-1)/2} (n-i) \binom{n}{i} \left(\frac{1}{2}\right)^n + \sum_{i=\frac{n-1}{2}+1}^{n} i \binom{n}{i} \left(\frac{1}{2}\right)^n \quad (21)$$

So, we continue with splitting the cases where n is even and n is odd. When n is an even number we find with some algebraic manipulation:

$$\begin{aligned}
E(X|N = 2k) &= \sum_{i=0}^{2k} \max(i, 2k-i) \binom{2k}{i} \left(\frac{1}{2}\right)^{2k} \\
&= \sum_{i=0}^{k} (2k-i) \binom{2k}{i} \left(\frac{1}{2}\right)^{2k} + \sum_{i=k+1}^{2k} i \binom{2k}{i} \left(\frac{1}{2}\right)^{2k} \\
&= \left(\frac{1}{2}\right)^{2k} 2k \left[\sum_{i=0}^{k} \binom{2k-1}{i} + \sum_{i=k}^{2k-1} \binom{2k-1}{i}\right] \\
&= \left(\frac{1}{2}\right)^{2k} 2k \left[2^{2k-1} + \binom{2k-1}{k}\right] \\
&= k + k \binom{2k-1}{k} \left(\frac{1}{2}\right)^{2k-1} \\
&= k + k \binom{2k}{k} \left(\frac{1}{2}\right)^{2k}
\end{aligned} \quad (22)$$

Likewise, for the odd case we find:

$$\begin{aligned}
E(X|N &= 2k+1) \\
&= \left(\frac{2k+1}{2}\right) + \left(\frac{2k+1}{2}\right) \binom{2k}{k} \left(\frac{1}{2}\right)^{2k}
\end{aligned} \quad (23)$$

With this result we can calculate EX:

$$\begin{aligned}
EX &= \sum_{n=0}^{\infty} E(X|N = n) P(N = n) \\
&= \sum_{k=0}^{\infty} [E(X|N = 2k) P(N = 2k) \\
&\quad + E(X|N = 2k+1) P(N = 2k+1)]
\end{aligned} \quad (24)$$

But first we must calculate $P(N = n)$, the probability that during a run of empty center sites a total of n particles were dropped in the border vertices. This is given by

$$P(N = n) = \frac{1}{3}\left(\frac{2}{3}\right)^n \quad (25)$$

which may be regarded trivial but can easily be calculated as follows. Let $T$ denote the time between two droppings in the center. Then we have

$$\begin{aligned}
P(N = n) &= \int_0^\infty P(N_T = n|T = t) P(T = t) dt \\
&= \int_0^\infty \frac{(2t)^n}{n!} e^{-2t} e^{-t} dt \\
&= \frac{2^n}{3^{n+1}} \int_0^\infty \frac{3^{n+1}}{n!} t^{(n+1)-1} e^{-3t} dt \\
&= \frac{1}{3}\left(\frac{2}{3}\right)^n
\end{aligned} \quad (26)$$

Continuing with Equation (24) we can now write:

$$EX = \sum_{n=0}^{\infty} E(X|N = n) P(N = n) \quad (27)$$

$$= \sum_{k=0}^{\infty} \left[ \left( k + k \binom{2k}{k} \left(\frac{1}{2}\right)^{2k} \right) \frac{1}{3} \left(\frac{2}{3}\right)^{2k} \right.$$
$$\left. + \left( \left(\frac{2k+1}{2}\right) + \left(\frac{2k+1}{2}\right) \binom{2k}{k} \left(\frac{1}{2}\right)^{2k} \right) \frac{1}{3} \left(\frac{2}{3}\right)^{2k+1} \right]$$
$$= \sum_{k=0}^{\infty} \left[ \left( k + \frac{2}{3} \frac{2k+1}{2} \right) \frac{1}{3} \left(\frac{2}{3}\right)^{2k} \right.$$
$$\left. + \frac{1}{3} k \binom{2k}{k} \left(\frac{1}{3}\right)^{2k} + \frac{2}{9} \left(\frac{2k+1}{2}\right) \binom{2k}{k} \left(\frac{1}{3}\right)^{2k} \right]$$
$$= \frac{1}{9} \sum_{k=0}^{\infty} \left(\frac{4}{9}\right)^{k} + \frac{5}{9} \sum_{k=0}^{\infty} k \left(\frac{4}{9}\right)^{k}$$
$$+ \frac{1}{3} \sum_{k=0}^{\infty} k \binom{2k}{k} \left(\frac{1}{9}\right)^{k} + \frac{2}{9} \sum_{k=0}^{\infty} k \binom{2k}{k} \left(\frac{1}{9}\right)^{k}$$
$$+ \frac{1}{9} \sum_{k=0}^{\infty} \binom{2k}{k} \left(\frac{1}{9}\right)^{k}$$
$$= \frac{1}{5} + \frac{4}{5} + \frac{5}{9} \frac{2/9}{(5/9)^{3/2}} + \frac{1}{9} \frac{3}{\sqrt{5}} = 1 + \frac{1}{\sqrt{5}}$$

where we used the identities (with $x = 1/9$)

$$\sum_{k=0}^{\infty} \binom{2k}{k} x^{k} = \frac{1}{\sqrt{1-4x}} \quad (28)$$

and

$$\sum_{k=0}^{\infty} k \binom{2k}{k} x^{k} = \frac{2x}{\sqrt{(1-4x)^3}} \quad (29)$$

where (28) can be easily checked by writing down the Taylor series expansion for $(1-4x)^{-1/2}$ while (29) follows immediately from (28) by differentiation with respect to x.

The first part of Theorem 2 follows now from Formula (17) and Lemma 2.

$$\lim_{r \to \infty} \rho_{\infty}(0, r) = \frac{1}{1 + EX} = \frac{1}{2 + \frac{1}{5}\sqrt{5}} = \frac{10 - \sqrt{5}}{19} \quad (30)$$

The proof of the second part of the theorem, i.e. that the limiting densities for the position left and right to the center are identical to the center, goes as follows. The expected total number of arrivals on one of the borders is equal to the number of arrivals in the center. This means that on average between two arrivals at position 0 there will be an arrival on the border site too. The distance between two particles which arrive at the center is on average EX. At the border site we expect one particle too. Therefore, the expected end-density on a border site must be $1/(1 + EX)$. This is indeed the same density as in the center.